\documentclass[a4paper,10pt,reqno]{amsart}
\usepackage[utf8]{inputenc}
\usepackage{graphicx}
\usepackage{amsfonts,amsmath,amssymb,amsthm}

\newcommand{\rmd}{\mathrm{d}}
\newcommand{\rme}{\mathrm{e}}
\newcommand{\rmO}{\mathrm{O}}
\newcommand{\Rset}{\mathbb{R}}

\newtheorem{lem}{Lemma}
\newtheorem{corol}{Corollary}
\newtheorem{defn}{Definition}
\newtheorem{thm}{Theorem}

\begin{document}

\title{Towards non-linear quadrature formulae}
\author[G.~M.~von~Hippel]{Georg M. von Hippel}
\address{Institut für Kernphysik \& PRISMA\textsuperscript{+} Cluster of
Excellence \\ University of Mainz \\ 55099 Mainz \\ Germany }



\date{12th May 2026}

\begin{abstract}
Prompted by an observation about the integral of exponential functions of the
form $f(x)=\lambda\rme^{\alpha x}$, we investigate the possibility to exactly
integrate families of functions generated from a given function by scaling or
by affine transformations of the argument using nonlinear generalizations of
quadrature formulae. The main result of this paper is that such formulae can be
explicitly constructed for a wide class of functions, and have the same
accuracy as Newton-Cotes formulae based on the same nodes, with the latter
emerging as the linear case of our general formalism. We also derive explicit
bounds on the error of the nonlinear quadrature formulae, which in the linear
case devolve into the well-known bounds for Newton-Cotes formulae.
\end{abstract}
\maketitle

\section{Introduction}

One of the most basic tasks in numerical analysis is the approximate evaluation
of definite integrals
\begin{equation}
I[f] = \int_a^{a+h}f(x)\,\rmd x
\end{equation}
by quadrature formulae
\begin{equation}
\hat{I}[f] = h q(\hat{f})
\label{eq:quadrature}
\end{equation}
where $\hat{f}=(f(a+\xi_0 h),\ldots,f(a+\xi_{n-1}h))\in\Rset^n$ and $q:\Rset^n\to\Rset$
such that
\begin{equation}
\left|\hat{I}[f]-I[f]\right|\le C_f h^r
\end{equation}
for some $r>1$ and some $f$-dependent constant $C_f$.
Since one of the most fundamental
properties of integration is its linearity, typical quadrature formulae consist
of taking linear combinations $q(\hat{f})=\langle\omega,\hat{f}\rangle$
of the values of the integrand at specific values
of the integration variable. Depending on whether only the linear coefficients
(weights) $\omega_i$ or also the abscissae (nodes) $\xi_i$ are adjusted to minimize
the error made in the numerical evaluation, one gets families of quadrature formulae
such as the Newton-Cotes or Gaussian quadrature formulae. These are then exact on
polynomials of a given degree.

However, in many applications the case arises that one needs to evaluate the integral of
some function given only in terms of sampled values, usually at regularly
spaced points, which is known to be very close to a family of
functions (other than polynomials) whose integrals are known analytically. As
a example, consider a function $f$ which is known to be very close to an
exponential. This arises for instance in the analysis of time series known to
consist of a sum of exponentially decaying components using the Padé-Laplace
method \cite{Yeramian:1987nature}, where the parameters $E_n$, $A_n$ of a
function $f$ known to be of the form
\[
f(t) = \sum_{n=1}^N A_n\rme^{-E_n t}
\]
but given only in terms of a finite number of measured values $f(x_i)$,
$x_i=ih$, $i<M$, $M>N$ are obtained as
the poles and residues of a Padé approximant to the Laplace transform of $f$,
i.e. the power series with coefficients
\begin{equation}
\left.\frac{\rmd^k}{\rmd p^k}\mathcal{L}[C](p)\right|_{p=p_0} = \int_0^\infty
\rmd t\,(-t)^k C(t)\rme^{-p_0t}.
\end{equation}
Another example arises in certain calculations in theoretical high-energy
physics \cite{Bernecker:2011gh}, such as the determination of the so-called
anomalous magnetic moment through an integral
\[
a^{\mathrm{HVP}}_\mu = \left(\frac{\alpha}{\pi}\right)^2 \int_0^\infty K(t)
G(t)\,dt
\]
where $\alpha$ is an experimentally measured physical constant and $K$ is an
analytically-known kernel function, but while the function $G(t)$ is known on
general grounds to be of the form
\begin{equation}
G(t) = \sum_{n=1}^\infty A_n \rme^{-E_n t}
\end{equation}
it can only be determined by simulations yielding only its values $G(t_i)$ at
fixed times $t_i=ia<T$. In these cases, it becomes important to be able to get
a good handle on the numerical evaluation of integrals of exponentially
decaying functions known only by tabulated values at regularly spaced points.

The integral of a function $f$ satisfying $|f(x)-\lambda\rme^{\alpha x}|<\epsilon$
can be approximated by
\begin{equation}
\int_a^b f(x)\,\rmd x =
\frac{\lambda}{\alpha}\left(\rme^{\alpha b}-\rme^{\alpha a}\right)
  +\rmO\left(\epsilon(b-a)\right).
\end{equation}
This approximation is still useful even if $\lambda$ and $\alpha$ are not known
beforehand, because they can be estimated from $f$. Indeed, within the bound
given by $\epsilon$, we can replace
$\lambda\rme^{\alpha a}$ and $\lambda\rme^{\alpha b}$ by $f(a)$ and $f(b)$,
respectively, 
and estimate $\alpha$ from the numerical derivative of the logarithm of $f$ as
$\alpha = \frac{1}{h}\left(\log f(x+h)-\log f(x)\right) + \rmO(h)$. Putting
these ingredients together, we arrive at a non-linear quadrature formula
\begin{equation}
\int_a^b f(x)\,\rmd x \approx \frac{f(b)-f(a)}{\log f(b)-\log f(a)} (b-a).
\label{eq:exemplary}
\end{equation}
Such formulae are used on an \emph{ad hoc} basis by practitioners in various
scientific disciplines, but a theoretical basis beyond the heuristic
considerations sketched above is currently lacking. Our purpose here is to put
such formulae on a firmer mathematical footing by developing a theory of
non-linear quadrature formulae encompassing eq.~(\ref{eq:exemplary}) as its
exemplary case.

The main questions we aim to answer are
\begin{enumerate}
\item What can we say about the accuracy of non-linear quadrature formulae of
the type of eq.~(\ref{eq:exemplary})?
\item Is there a systematic way to achieve improved accuracy by using more than
two evaluations of $f$ also in the non-linear setting?
\end{enumerate}
Our main results can be summarized as
\begin{enumerate}
\item non-linear quadrature formulae that are exact on families of functions of
the form $f_{\alpha\beta}(x)=f^*(\alpha x+\beta)$ for some function $f^*$ have
accuracy $r=3$, i.e. comparable to the trapezoidal rule, which becomes the
better the more the integrand resembles the functions they are exact on,
\item under suitable conditions on the two-point non-linear quadrature rule,
a form of Romberg improvement can be performed on it to obtained a three-point
non-linear quadrature rule of accuracy $r=5$, i.e. comparable to Simpson's
rule.
\end{enumerate}

Even though the results we will obtain do not require more than standard
undergraduate analysis and therefore ought to be
well-known, there appears to be hardly any literature on the topic of
non-linear quadrature formulae, apart from two papers by Werner \cite{Werner:1978CMwA}
and Wuytack \cite{Wuytack:1975JCAM}, which relate to the use of rational interpolants
or Padé approximants instead of interpolating polynomials to integrate functions known
to have a singularity at one end of the integration interval. Our approach here
will be much more general.

We note that the formula of eq.~(\ref{eq:exemplary}) is exact for functions of
the form $f(x)=\lambda\rme^{\alpha x}$, which form a non-linear family that can
equivalently be expressed as $f(x)=\rme^{\alpha x+\beta}$. Moreover, swapping
the values of $f(a)$ and $f(a+h)$ does not change the value of the
approximation, and an overall factor can be pulled out of the approximation
since it cancels within the denominator.

Based on these observations, we will consider non-linear approximations to $I[f]$
by considering non-linear functions $q:\Rset^n\to\Rset$ in
eq.~(\ref{eq:quadrature}).
We will assume throughout that $q$ and $f$ are sufficiently smooth.
To make such approximations useful, one typically has to require that they
become exact in some limiting case. We therefore define several properties that
will become useful in the following:

\begin{defn}
A (non-linear) quadrature formula $\hat{I}$ is
\begin{enumerate}
\item \emph{exact} on a function $f^*$ if $\hat{I}[f^*]=I[f^*]$ for all $h>0$,
\item \emph{scalably exact} on a function $f^*$ if there exist
$\lambda_-<1<\lambda_+$ such that $\hat{I}[\lambda f^*]=I[\lambda f^*]$
for all $\lambda\in(\lambda_-;\lambda_+)$ and all $h>0$,
\item \emph{affinely exact} on a function $f^*$ if with
$f^*_{\alpha,\beta}(x)=f^*(\alpha x+\beta)$
we have $\hat{I}[f^*_{\alpha,\beta}]=I[f^*_{\alpha,\beta}]$
for all $\alpha,\beta\in\Rset$ and $h>0$,
\item \emph{symmetric} if $\xi_{n-1-k}=1-\xi_k$ for all $k\in\{0,\ldots,n-1\}$
and $q(f_{n-1},\ldots,f_0)=q(f_0,\ldots,f_{n-1})$ for all $f\in\Rset^n$, and
\item \emph{quasilinear} if $q(\lambda\hat{f})=\lambda q(\hat{f})$ for all
$\lambda\in\Rset$ and all $\hat{f}\in\Rset^n$.
\end{enumerate}
\end{defn}

The idea behind these definitions is that we will take exactness on a family of
target functions as our guide
as to the goodness of a quadrature rule (while noting that this has recently
been pointed out by Trefethen \cite{Trefethen:20210109501} to not be an entirely
reliable heuristic in the case of traditional linear quadrature rules), and we
will attempt to preserve at least some of the linear properties of integration
when acting on these target functions.

The main result of this paper is the following

\begin{thm}
\label{thm:errorbounds}
Let $n=2$, $\xi_0=0$, $\xi_1=1$, and let $\hat{I}$ be affinely
exact on some function $f^*$. Then
\begin{enumerate}
\item if $f^*\in C^3(\Rset)$ with $f^*(x)\not=0$ and
$f^{*}{}'(x)\not=0$ for all $x\in\Rset$, we have for $f\in C^3(\Rset)$
\begin{equation}
\label{eq:errbound1}
|\hat{I}[f]-I[f]|\le\frac{h^3}{6}\sup_{\xi\in[a;a+h]}\left|N[f](\xi)\right|
\end{equation}
where
\begin{align}
N[f](\xi) &= f''(\xi)[1-3 q^{(0,1)}(f(a),f(\xi))]
-3 f'(\xi)^2 q^{(0,2)}(f(a),f(\xi))\\
&-(\xi-a) \left[ f^{(3)}(\xi) q^{(0,1)}(f(a),f(\xi))
         +3 f'(\xi)f''(\xi)q^{(0,2)}(f(a),f(\xi))\right.\nonumber\\
         &~~~~~~~~\left.+f'(\xi)^3 q^{(0,3)}(f(a),f(\xi))\right]\nonumber
\end{align}
is a nonlinear function of $f$ satisfying $N[f^*_{\alpha,\beta}]=0$ for all
$\alpha,\beta\in\Rset$.
\item if $f^*\in C^2(\Rset\to R)$, $R\subseteq\Rset$, is bijective
with inverse function $f^*{}^{(-1)}$, whe have for $f\in C^2(\Rset \to R)$
\begin{equation}
\label{eq:errbound2}
\left|\hat{I}[f]-I[f]\right|\le\frac{h^3}{12}\sup_{\xi\in[a;a+h]}\left|f''(\xi)-\alpha_f^2f^*{}''(\alpha_f\xi+\beta_f)\right|
\end{equation}
where
\begin{align}
\alpha_f &= \frac{{f^*}^{(-1)}(f(a+h))-{f^*}^{(-1)}(f(a))}{h} \\
\beta_f &= \frac{(a+h) {f^*}^{(-1)}(f(a))-a {f^*}^{(-1)}(f(a+h))}{h}
\end{align}
\item if $f^*\in C^2(\Rset\to R)$, $R\subseteq\Rset$, is bijective
with inverse function $f^*{}^{(-1)}$ and satisfies
$Lf^*=0$, where $(Lu)(x)=-(p(x)u'(x))'+q(x)u(x)$
for $p(x),q(x)>0$ on $[a,a+h]$, we have for $f\in C^2(\Rset \to R)$
\begin{equation}
\label{eq:errbound3}
\left|\hat{I}[f]-I[f]\right|\le C_f
\sup_{\xi\in[a;a+h]}\left|L_{\alpha_f,\beta_f}f\right|
\end{equation}
where
\begin{equation}
(L_{\alpha,\beta}u)(x) = -(p(\alpha x+\beta)u'(x))'+\alpha^2 q(\alpha x+\beta)u(x)
\end{equation}
and $C_f=\int_a^{a+h}\int_a^{a+h} G_f(x,y)\,\rmd x\rmd y$ with $G_f(x,y)$ being the Green
function of $L_{\alpha_f,\beta_f}$ subject to Dirichlet boundary conditions on $[a;a+h]$.
\end{enumerate}
\end{thm}

as well as the explicit construction of the following

\begin{thm}
\label{thm:construct}
Let $f^*:\Rset\to R\subseteq\Rset$ be bijective with inverse function
$f^*{}^{(-1)}$, and let $F^*$ be an antiderivative of $f^*$. Define
\begin{equation}\label{eq:construct}
q_1(f_0,f_1) =
\frac{F^*(f^*{}^{(-1)}(f_1))-F^*(f^*{}^{(-1)}(f_0))}{f^*{}^{(-1)}(f_1)-f^*{}^{(-1)}(f_0)}.
\end{equation}
Then
\begin{equation}
\hat{I}_1[f]=h q_1(f(a),f(a+h))
\end{equation}
is symmetric and affinely exact on $f^*$, and with
\begin{equation}\label{eq:constructh4}
q_2(f_0,f_1,f_2) =
\frac{2}{3}\left(q_1(f_0,f_1))+q_1(f_1,f_2)\right)-\frac{1}{3}q_1(f_0,f_2)
\end{equation}
the three-point non-linear quadrature formula with $\xi_0=0$,
$\xi_1=\frac{1}{2}$, $\xi_2=1$ given by
\begin{equation}
\hat{I}_2[f] = h q_2(f(a),f(a+\frac{h}{2}),f(a+h))
\end{equation}
is symmetric and affinely exact on $f^*$ and satisfies
\begin{equation}
\left|\hat{I}_2[f]-I[f]\right|\le \frac{h^5}{2880}
\sup_{\xi\in[a;a+h]}|f^{(4)}(\xi)-\Omega_f|
\end{equation}
where
\begin{align}
\Omega_f&=\frac{480}{h^4}\Bigg(f(a)+4f(a+\frac{h}{2})+f(a+h)+2q_1(f(a),f(a+h))\\
&~~~~-4\left(q_1(f(a),f(a+\frac{h}{2}))+q_1(a+\frac{h}{2},f(a+h))\right)\Bigg)\nonumber
\end{align}
remains finite as $h\to0$.
\end{thm}

As an excursion, we will consider the traditional Newton-Cotes formulae as
special (linear) cases of the general (non-linear) case and rederive some
well-known results in this way. Finally, we give some explicit examples and perform
some numerical experiments to investigate the potential usefulness and
limitations of non-linear quadrature rules.

\section{Two-point non-linear quadrature rules}

First, we show that scalably exact quadratures have at least no worse order
than the trapezoidal rule on general functions:
\begin{lem}\label{lem:noworse}
Let $n=2$, $\xi_0=0$, $\xi_1=1$.
If there exists a function $f^*\in C^2(a,a+h)$ with $f^*(a)\not=0$ and $f^{*}{}'(a)\not=0$
such that $\hat{I}$ is scalably exact on $f^*$,
then $\left|\hat{I}[f]-I[f]\right|=o(h^2)$ for all $f\in C^2(a,a+h)$.
\begin{proof}
The Taylor expansion in $h$ of the exact integration is given by
\begin{equation}
I[f] = f(a) h + \frac{1}{2}f'(a)h^2 + \frac{1}{6}f''(a)h^3 + o(h^3)
\end{equation}
whereas that of the two-point quadrature formula is given by
\begin{equation}
\hat{I}[f] = q(\bar{f}) h + q^{(0,1)}(\bar{f}) f'(a) h^2 +
\frac{1}{2}\left(q^{(0,2)}(\bar{f})[f'(a)]^2+q^{(0,1)}(\bar{f})f''(a)\right)h^3 + o(h^3)
\end{equation}
where $\bar{f}=(f(a),f(a))$ and we understand $q(\bar{f})$ to denote
$\lim_{f_b\to f(a)} q(f(a),f_b)$ in the case where $q(\bar{f})$ itself is
ill-defined. In order for these to be identical for all
$f=\lambda f^*$, we need to have $q(\bar{f})=f(a)$ and
$q^{(0,1)}(\bar{f})=\frac{1}{2}$ for all $\bar{f}$.
Hence, we have for arbitrary $f\in C^2(a,a+h)$ that
\begin{equation}
\left|\hat{I}[f]-I[f]\right| =
\left|\frac{1}{2}q^{(0,2)}(\bar{f})[f'(a)]^2+\frac{1}{12}f''(a)\right|h^3 + o(h^3) =
o(h^2).
\end{equation}
\end{proof}
\end{lem}

Noting that the exactness of $\hat{I}$ on functions of the form $\lambda f^*$
requires linearity of $q$ in the vicinity of $\bar{f}^*$ since $\hat{I}[\lambda
f^*]=I[\lambda f^*]=\lambda I[f^*]=\lambda\hat{I}[f^*]$, we find relationships
between the partial derivatives of $q$:
\begin{lem}
Let $\hat{I}$, $f^*$ be as in Lemma \ref{lem:noworse}. Then
\begin{align}
q^{(1,0)}(\bar{f}^*) + q^{(0,1)}(\bar{f}^*) &= 1, \nonumber\\
q^{(0,2)}(\bar{f}^*) + q^{(1,1)}(\bar{f}^*) &= 0, \\
q^{(0,2)}(\bar{f}^*) + 2q^{(1,1)}(\bar{f}^*) + q^{(2,0)}(\bar{f}^*) &= 0
   \nonumber\\
q^{(0,3)}(\bar{f}^*)
+3q^{(1,2)}(\bar{f}^*) + 3q^{(2,1)}(\bar{f}^*)+q^{(3,0)}(\bar{f}^*) &= 0,
   \nonumber\\
q^{(0,3)}(\bar{f}^*) + 2q^{(1,2)}(\bar{f}^*)+q^{(2,1)}(\bar{f}^*) &= 0.
   \nonumber
\end{align}
\begin{proof}
We expand the exact integral (which is of course linear) in powers of $h$
\begin{equation}
I[\lambda f^*] = h\lambda f^*(a) + \frac{1}{2}h^2\lambda
f^*{}'(a)+\frac{1}{6}h^3\lambda f^*{}''(a)+o(h^3)
\end{equation}
and perform a double expansion of $\hat{I}[\lambda f^*]$ into powers of $h$ and
$\lambda-1$,
\begin{align}
\hat{I}[\lambda f^*] &=
h\Big[q(\bar{f}^*)+(\lambda-1)f(a)\left(q^{(1,0)}(\bar{f}^*) +
q^{(0,1)}(\bar{f}^*)\right) \\
&\qquad\qquad +\frac{(\lambda-1)^2}{2}f(a)^2\left(q^{(0,2)}(\bar{f}^*) +
2q^{(1,1)}(\bar{f}^*) + q^{(2,0)}(\bar{f}^*)\right) \nonumber\\
&\qquad +\frac{(\lambda-1)^3}{6}f(a)^3\left(q^{(0,3)}(\bar{f}^*)
+3q^{(1,2)}(\bar{f}^*) + 3q^{(2,1)}(\bar{f}^*)+q^{(3,0)}(\bar{f}^*)\right)\Big] \nonumber\\
&+\frac{1}{2}h^2 f^*{}'(a)\Big[\lambda +
2\lambda(\lambda-1)f(a)\left(q^{(0,2)}(\bar{f}^*) +
q^{(1,1)}(\bar{f}^*)\right)\nonumber\\
&\qquad\qquad +(\lambda-1)^2 \lambda f(a)^2\left(q^{(0,3)}(\bar{f}^*) +
2q^{(1,2)}(\bar{f}^*)+q^{(2,1)}(\bar{f}^*)\right)\Big]\nonumber\\
&+o(h^2)+o(h(\lambda-1)^3)+o(h^2(\lambda-1)^2)\nonumber
\end{align}
and equality for all $\lambda$ and $h$ is only possible if the given relations
hold.
\end{proof}
\end{lem}

We therefore find that scalably exact non-linear quadratures locally resemble the
trapezoidal rule:

\begin{corol} Let $n=2$, $\xi_0=0$, $\xi_1=1$, and let $\hat{I}$ be scalably
exact on $f^*\in C^2(a,a+h)$ with $f^*{}'(a)\not=0$.
Then
\begin{equation}
q(f(a),f(b)) = \frac{f(a)+f(b)}{2} - \frac{1}{12} \frac{f^*{}''(a)}{[f^*{}'(a)]^2}
(f(b)-f(a))^2 + o\left(|f(b)-f(a)|^2\right).
\end{equation}
\begin{proof}
Since $\hat{I}[f^*]=I[f^*]$, we must have
\begin{equation}
\frac{1}{2}\left(q^{(0,2)}(\bar{f}^*)[f^{*'}(a)]^2+q^{(0,1)}(\bar{f}^*)f^*{}''(a)\right)
= \frac{1}{6}f^*{}''(a)
\end{equation}
and using $q^{(0,1)}(\bar{f})=\frac{1}{2}$ then yields
\begin{equation}
q^{(0,2)}(\bar{f}^*) = -\frac{1}{6} \frac{f^*{}''(a)}{[f^*{}'(a)]^2}.
\end{equation}
Substituting this into the Taylor expansion of $q(f(a),f(b))$ around
$f(b)=f(a)$ and using the relations between the partial derivatives of $q$
found above yields the result.
\end{proof}
\end{corol}

Essentially identical results can be shown for affinely exact non-linear
quadrature formulae:
\begin{lem}\label{lem:aff}
Let $n=2$, $\xi_0=0$, $\xi_1=1$.
If there exists a function $f^*\in C^2(\Rset)$ with $f^*(x)\not=0$ and
$f^{*'}(x)\not=0$ for all $x\in\Rset$
such that $\hat{I}$ is affinely exact on $f^*$, 
then $\left|\hat{I}[f]-I[f]\right|=o(h^2)$ for all $f\in C^2(a,a+h)$. Moreover,
\begin{align}
q^{(1,0)}(\bar{f}^*) + q^{(0,1)}(\bar{f}^*) &= 1, \nonumber\\
q^{(0,2)}(\bar{f}^*) + q^{(1,1)}(\bar{f}^*) &= 0, \\
q^{(0,2)}(\bar{f}^*) + 2q^{(1,1)}(\bar{f}^*) + q^{(2,0)}(\bar{f}^*) &= 0,
   \nonumber
\end{align}
and for $f\in C^3(a,a+h)$ also
\begin{align}
q^{(0,3)}(\bar{f}^*)
+3q^{(1,2)}(\bar{f}^*) + 3q^{(2,1)}(\bar{f}^*)+q^{(3,0)}(\bar{f}^*) &= 0,
   \\
q^{(0,3)}(\bar{f}^*) + 2q^{(1,2)}(\bar{f}^*)+q^{(2,1)}(\bar{f}^*) &= 0.
   \nonumber
\end{align}
\begin{proof}
We start by noting that the Taylor expansion of
$I[f^*_{\alpha,\beta}]-\hat{I}[f^*_{\alpha,\beta}]$ in $h$ is
\begin{equation}
I[f^*_{\alpha,\beta}]-\hat{I}[f^*_{\alpha,\beta}] = h \left(f^*(\alpha a+\beta)
-q(\bar{f}^*_{\alpha,\beta})\right) +
\frac{\alpha h^2}{2}f^{*}{}'(\alpha a+\beta)
\left(1-2q^{(0,1)}(\bar{f}^*_{\alpha,\beta})\right)+ o(h^2)
\end{equation}
which can only vanish for all $h,\alpha,\beta$ if for all values of
$\bar{f}$ with $f(a)$ in the range of $f^*$ the equalities $q(\bar{f})=f(a)$
and $q^{(0,1)}(\bar{f})=\frac{1}{2}$ hold, implying
$\left|\hat{I}[f]-I[f]\right|=o(h^2)$ as in Lemma \ref{lem:noworse}.
Furthermore, the Taylor expansions around $\alpha=0$ of the coefficients of $h$
and $h^2$ are 
\begin{align}
f^*(\alpha a+\beta)-q(\bar{f}^*_{\alpha,\beta}) 
&= f^*(\beta)-q(\bar{f}^*_{0,\beta})+ \alpha a
f^*{}'(\beta)\left(1-q^{(0,1)}(\bar{f}^*_{0,\beta})-q^{(1,0)}(\bar{f}^*_{0,\beta})\right)
\\
&+\frac{(\alpha a)^2}{2}
\Big[f^*{}''(\beta)\left(1-q^{(0,1)}(\bar{f}^*_{0,\beta})-q^{(1,0)}(\bar{f}^*_{0,\beta})\right)\nonumber\\
&~~~~~ - [f^*{}'(\beta)]^2\left(q^{(0,2)}(\bar{f}^*_{0,\beta}) +
2q^{(1,1)}(\bar{f}^*_{0,\beta}) +
q^{(2,0)}(\bar{f}^*_{0,\beta})\right)\Big] \nonumber\\
&+\frac{(\alpha a)^3}{6}
\Big[f^*{}^{(3)}(\beta)\left(1-q^{(0,1)}(\bar{f}^*_{0,\beta})-q^{(1,0)}(\bar{f}^*_{0,\beta})\right)\nonumber\\
&~~~~~ - 3f^*{}'(\beta)f^*{}''(\beta)\left(q^{(0,2)}(\bar{f}^*_{0,\beta}) +
2q^{(1,1)}(\bar{f}^*_{0,\beta}) +
q^{(2,0)}(\bar{f}^*_{0,\beta})\right)\nonumber\\
&~~~~~ -[f^*{}'(\beta)]^3\left(q^{(0,3)}(\bar{f}^*_{0,\beta})+3q^{(1,2)}(\bar{f}^*_{0,\beta})+3q^{(2,1)}(\bar{f}^*_{0,\beta})+q^{(3,0)}(\bar{f}^*_{0,\beta})\right)
\Big]\nonumber \\
&+ o(\alpha^3)\nonumber
\end{align}
and
\begin{align}
\left(1-2q^{(0,1)}(\bar{f}^*_{\alpha,\beta})\right) f^*{}'(\alpha a+\beta) &= 
f^*{}'(\beta)q^{(0,1)}(\bar{f}^*_{0,\beta})\\
&+\alpha a \Big[f^*{}''(\beta)q^{(0,1)}(\bar{f}^*_{0,\beta}) \nonumber\\
&~~~~~ +[f^*{}'(\beta)]^2\left(q^{(0,2)}(\bar{f}^*_{0,\beta})+q^{(1,1)}(\bar{f}^*_{0,\beta})\right)\Big]
\nonumber \\
&+\frac{(\alpha a)^2}{2} \Big[
f^*{}^{(3)}(\beta)q^{(0,1)}(\bar{f}^*_{0,\beta})\nonumber \\
&~~~~~ +3 f^*{}'(\beta)f^*{}''(\beta)\left(q^{(0,2)}(\bar{f}^*_{0,\beta}) +
q^{(1,1)}(\bar{f}^*_{0,\beta})\right)\nonumber\\
&~~~~~ +[f^*{}'(\beta)]^3\left(q^{(0,3)}(\bar{f}^*_{0,\beta})+2q^{(1,2)}(\bar{f}^*_{0,\beta})+q^{(2,1)}(\bar{f}^*_{0,\beta})\right)
\Big]\nonumber \\
&+ o(\alpha^2),\nonumber
\end{align}
which must vanish identically, whence the relations immediately follow.
\end{proof}
\end{lem}

The equations for the partial derivatives of scalably or affinely exact $q$ are
readily solved, leading to the following

\begin{corol}\label{cor:pders}
Let $n=2$, $\xi_0=0$, $\xi_1=1$, and let $\hat{I}$ be scalably or affinely
exact on some function $f^*\in C^2(\Rset)$ with $f^*(a)\not=0$ and
$f^{*'}(a)\not=0$. Then
\begin{align}
q(\bar{f}) &= f(a), \nonumber\\
q^{(1,0)}(\bar{f}) &= q^{(0,1)}(\bar{f}) = \frac{1}{2}, \\
q^{(2,0)}(\bar{f}) &= -q^{(1,1)}(\bar{f}) = q^{(0,2)}(\bar{f}). \nonumber
\end{align}
If $\hat{I}$ is moreover symmetric and $f^*\in C^3(\Rset)$, then also
\begin{equation}
q^{(3,0)}(\bar{f}) = -3 q^{(1,2)}(\bar{f})
 = -3 q^{(2,1)}(\bar{f}) = q^{(0,3)}(\bar{f}).
\end{equation}
\qed
\end{corol}

In the following, we will find it useful to have shown the following

\begin{lem}
\label{lem:greenbound}
For $f\in C^r(\Rset)$, let $p_f\in C^r(\Rset)$ be defined such that
$p_f(a+\xi_k h)=f(a+\xi_k h)$, $k=0,\ldots,n-1$, $p_f'(a+\xi_k h)=f'(a+\xi_k
h)$, $k\in H\subset\{1,\ldots,n-2\}$, and $\hat{I}[p_f]=I[p_f]$.
Let $D$ be a differential operator of
order $r\ge n$ such that $Du=0$ on $[a;a+h]$ with boundary conditions
$u(a+\xi_k h)=0$, $k=0,\ldots,n-1$, $u'(a+\xi_k h)=0$,
$k\in H\subset\{1,\ldots,n-2\}$, $|H|=r-n$, has no non-trivial solutions.
Then we have
\begin{equation}
|\hat{I}[f]-I[f]|\le C \sup_{\xi\in[a;a+h]} |Df(\xi)-Dp_f(\xi)|
\end{equation}
where
\begin{equation}
C = \int_a^{a+h} \int_a^{a+h} |G(x,y)|\,\rmd x\,\rmd y
\end{equation}
is the double integral of the absolute value of the Green function $G$ of $D$.
\begin{proof}
Note that by construction $\hat{I}[p_f]=\hat{I}[f]$, whence we have
$|\hat{I}[f]-I[f]|=|I[p_f]-I[f]|=|I[e_f]|\le I[|e_f|]$ for $e_f=f-p_f$. Now $De_f=Df-Dp_f$ by
definition, with solution $e_f(x) = \int_a^{a+h} G(x,y)(Df(y)-Dp_f(y))\,\rmd y$. We
therefore have
\begin{align}
I[|e_f|] &= \int_a^{a+h} \left|\int_a^{a+h}
G(x,y)(Df(y)-Dp_f(y))\,\rmd y\right|\,\rmd x\nonumber\\
&\le \int_a^{a+h} \int_a^{a+h}
|G(x,y)|\,|Df(y)-Dp_f(y)|\,\rmd x\,\rmd y
\end{align}
and the given bound follows.
\end{proof}
\end{lem}

We now proceed to prove the error bounds of Theorem \ref{thm:errorbounds}.
Let $n=2$, $\xi_0=0$, $\xi_1=1$, and let $\hat{I}$ be affinely
exact on some function $f^*$.

The most naive error bound that can be derived is given by
\begin{equation}
|\hat{I}[f]-I[f]|\le\frac{h^3}{6}\sup_{\xi\in[a;a+h]}\left|N[f](\xi)\right|
\end{equation}
for $f^*,f\in C^3(\Rset)$ with $f^*(x)\not=0$ and
$f^{*'}(x)\not=0$ for all $x\in\Rset$, where
\begin{align}
N[f](\xi) &= f''(\xi)[1-3 q^{(0,1)}(f(a),f(\xi))]
-3 f'(\xi)^2 q^{(0,2)}(f(a),f(\xi))\\
&-(\xi-a) \left[ f^{(3)}(\xi) q^{(0,1)}(f(a),f(\xi))
         +3 f'(\xi)f''(\xi)q^{(0,2)}(f(a),f(\xi))\right.\nonumber\\
         &\left.\qquad\qquad+f'(\xi)^3 q^{(0,3)}(f(a),f(\xi))\right]\nonumber
\end{align}
is a nonlinear function of $f$ satisfying $N[f^*_{\alpha,\beta}]=0$ for all
$\alpha,\beta\in\Rset$.
\begin{proof}
Taylor-expanding the error functional $\hat{I}[f]-I[f]$ in $h$
with the remainder term in Lagrange form yields
\begin{equation}
\hat{I}[f]-I[f] = h\left(q(\bar{f})-f(a)\right)
+\frac{h^2}{2} \left(2 q^{(0,1)}(\bar{f}) -1 \right)f'(a)
+\frac{h^3}{6} N[f](\xi)
\end{equation}
for some $\xi\in[a;a+h]$,
and for an affinely exact formula the first two terms vanish, leaving
\begin{equation}
|\hat{I}[f]-I[f]| = \frac{h^3}{6} |N[f](\xi)|
\end{equation}
whence the error bound immediately follows. For functions of the form
$f=f^*_{\alpha,\beta}$, the error functional vanishes exactly for all values
of $a$ and $h$, which means that $N[f^*_{\alpha,\beta}]$ must have a zero in
any given interval $[a;a+h]$, and therefore must vanish identically.
\end{proof}

An error bound that does not require third derivatives is given by
\begin{equation}
\left|\hat{I}[f]-I[f]\right|\le\frac{h^3}{12}\sup_{\xi\in[a;a+h]}\left|f''(\xi)-\alpha_f^2f^*{}''(\alpha_f\xi+\beta_f)\right|
\end{equation}
for bijective $f^*\in C^2(\Rset\to R)$, $R\subseteq\Rset$,
with inverse function $f^*{}^{(-1)}$, and $f\in C^2(\Rset \to R)$,
where
\begin{align}
\alpha_f &= \frac{{f^*}^{(-1)}(f(a+h))-{f^*}^{(-1)}(f(a))}{h} \\
\beta_f &= \frac{(a+h) {f^*}^{(-1)}(f(a))-a {f^*}^{(-1)}(f(a+h))}{h}
\end{align}
\begin{proof}
We note that $f^*(\alpha_fa+\beta_f)=f(a)$, $f^*(\alpha_f(a+h)+\beta_f)=f(a+h)$
by construction, and
$\hat{I}[f^*_{\alpha_f,\beta_f}]=I[f^*_{\alpha_f,\beta_f}]$ from the
affine exactness of $\hat{I}$. We can thus apply lemma \ref{lem:greenbound}
with $p_f=f^*_{\alpha_f,\beta_f}$, $r=2$, $D=\frac{\rmd^2}{\rmd x^2}$, whence
$Df(\xi)=f''(\xi)$, $Dp_f(\xi)=\alpha_f^2f^*{}''(\alpha_f\xi+\beta_f)$, and
$G(x,y)=\frac{(a-y) (a+h-x) \theta (x-y)+(a-x) (a+h-y)\theta (y-x)}{h}$, thus
$C=\frac{h^3}{12}$.
\end{proof}

Another error bound that looks more conventional (as in not explicitly
involving the difference of two functions) and requires only second derivatives
can be given under somewhat stronger requirements on $f^*$:
\begin{equation}
\left|\hat{I}[f]-I[f]\right|\le C_f
\sup_{\xi\in[a;a+h]}\left|L_{\alpha_f,\beta_f}f\right|
\end{equation}
for bijective $f^*\in C^2(\Rset\to R)$, $R\subseteq\Rset$,
which satisfies $Lf^*=0$, where $(Lu)(x)=-(p(x)u'(x))'+q(x)u(x)$
for $p(x),q(x)>0$ on $[a,a+h]$, and $f\in C^2(\Rset \to R)$, where
\begin{equation}
(L_{\alpha,\beta}u)(x) = -(p(\alpha x+\beta)u'(x))'+\alpha^2 q(\alpha x+\beta)u(x)
\end{equation}
and $C_f=\int_a^{a+h}\int_a^{a+h} |G_f(x,y)|\,\rmd x\rmd y$ with $G_f(x,y)$ being the Green
function of $L_{\alpha_f,\beta_f}$ subject to Dirichlet boundary conditions on $[a;a+h]$.
\begin{proof}
We first note that $L_{\alpha,\beta}f^*_{\alpha,\beta}=-(p(\alpha
x+\beta)\alpha f^*{}'(\alpha x+\beta))'+\alpha^2 q(\alpha x+\beta)f^*(\alpha
x+\beta) = \alpha^2[-(p(y)f^*{}'(y))'+q(y)f^*(y)]=\alpha^2(Lf^*)(y)=0$
with $y=\alpha x+\beta$. Applying lemma \ref{lem:greenbound} 
with $p_f=f^*_{\alpha_f,\beta_f}$, $r=2$, and $D=L_{\alpha_f,\beta_f}$ then
yields $Df=L_{\alpha_f,\beta_f}f$, $Dp_f=0$, and $C=C_f$.
\end{proof}
We note that $C_f=O(h^3)$ from its definition.
This completes the proof of Theorem \ref{thm:errorbounds}.

At least for a certain class of functions, symmetric affinely exact non-linear
quadrature formulae are readily constructed, which is the first part of Theorem
\ref{thm:construct}:

Let $f^*:\Rset\to R\subseteq\Rset$ be bijective with inverse function
$f^*{}^{(-1)}$, and let $F^*$ be an antiderivative of $f^*$. Define
\begin{equation}
q_1(f_0,f_1) =
\frac{F^*(f^*{}^{(-1)}(f_1))-F^*(f^*{}^{(-1)}(f_0))}{f^*{}^{(-1)}(f_1)-f^*{}^{(-1)}(f_0)}.
\end{equation}
Then
\begin{equation}
\hat{I}_1[f]=h q_1(f(a),f(a+h))
\end{equation}
is symmetric and affinely exact on $f^*$.

\begin{proof}
We have
\begin{equation}
q_1(f^*(\alpha a+\beta),f^*(\alpha (a+h)+\beta)) = \frac{F^*(\alpha
(a+h)+\beta)-F^*(\alpha a + \beta)}{\alpha(a+h+\beta)-(\alpha a+\beta)}
\end{equation}
and hence
\begin{equation}
\hat{I}[f^*_{\alpha,\beta}] = \frac{F^*(\alpha
(a+h)+\beta)-F^*(\alpha a + \beta)}{\alpha} = I[f^*_{\alpha,\beta}]
\end{equation}
as required. The symmetry of $q_1$ under an interchange of its two arguments is
readily apparent.
\end{proof}

\section{Towards higher-order non-linear quadrature rules}

In order to achieve higher order without moving the nodes $\xi_i$ so as to
require an evaluation of $f$ away from the equally spaced sampling points at
which it is known in typical applications, quadrature rules with additional
points have to be considered. Here we will constrain ourselves to the case of
three-point rules of the form $\hat{I}[f]=h q(f(a),f(a+h/2),f(a+h))$.

Given a suitable two-point non-linear quadrature rule, one can readily
construct a three-point non-linear quadrature rule of higher order using what
is essentially Romberg improvement:

Let $\hat{I}_1[f]=h q_1(f(a),f(a+h))$ be symmetric and affinely exact on a
bijective function $f^*\in C^4(\Rset\to R)$, $R\subseteq\Rset$. Then for
\begin{equation}
\label{eq:q4def}
q_2(f_0,f_1,f_2) =
\frac{2}{3}\left(q_1(f_0,f_1))+q_1(f_1,f_2)\right)-\frac{1}{3}q_1(f_0,f_2)
\end{equation}
the three-point non-linear quadrature formula with $\xi_0=0$,
$\xi_1=\frac{1}{2}$, $\xi_2=1$ given by
\begin{equation}
\hat{I}_2[f] = h q_2(f(a),f(a+\frac{h}{2}),f(a+h))
\end{equation}
is symmetric and affinely exact on $f^*$ with $o(h^4)$ errors.
\begin{proof}
Since $q_1$ is a symmetric affinely exact non-linear quadrature formula, its
partial derivatives at each order are given by Corollary \ref{cor:pders}.
We form the linear combination
\begin{align}
q_2(\hat{f}) &= \alpha_1 q_1(f(a),f(a+h))\\
&+\alpha_2\left(q_1(f(a),f(a+\frac{1}{2}h))+q_1(f(a+\frac{1}{2}h),f(a+h))\right)\nonumber
\end{align}
and determine the weights $\alpha_i$ from expanding
\begin{align}
 h q_2(\hat{f}) &=
   \left(\alpha _1+2 \alpha _2\right) h f(a)+\frac{1}{2} \left(\alpha _1+2 \alpha
   _2\right) h^2 f'(a)\\\nonumber&+\frac{1}{8} h^3 \Big(2 \left(2 \alpha _1+\alpha _2\right)
   f'(a)^2 q_1{}^{(0,2)}(f(a),f(a))\\\nonumber&\;\;\;\;\;+\left(2 \alpha _1+3 \alpha _2\right)
   f''(a)\Big)\\\nonumber&+\frac{1}{48} h^4 \Big(
   4 \left(2 \alpha _1+\alpha _2\right) f'(a)^3 q_1{}^{(0,3)}(f(a),f(a))
\\\nonumber&\;\;\;\;\;+12\left(2 \alpha _1+\alpha _2\right) f''(a)q_1{}^{(0,2)}(f(a),f(a))
   \\\nonumber&\;\;\;\;\;+ \left(4\alpha _1 +5 \alpha _2\right) f^{(3)}(a) \Big)
+o\left(h^4\right)
\nonumber
\end{align}
and equating this with
\begin{equation}
I[f] = h f(a)+\frac{1}{2} h^2 f'(a)+\frac{1}{6} h^3 f''(a)+\frac{1}{24} h^4
   f^{(3)}(a)+o\left(h^4\right)
\end{equation}
which yields the solution $\alpha_1=-\frac{1}{3}$, $\alpha_2=\frac{2}{3}$.

Since each of
the approximations across subintervals is affinely exact, so is their linear
combination. The symmetry of $q_2$ follows from that of $q_1$ by inspection.
\end{proof}

We note that eliminating the $O(h^3)$ term from the error also eliminated the
$O(h^4)$ term due to the symmetry of the two-point rule, in complete analogy to
what happens in the linear case \cite{vonPetersdorff:1993AMM}.

Deriving bounds on the error becomes more difficult in the three-point case,
since the condition of affine exactness gives us only two parameters to create
an interpolating function that is integrated exactly. Without strengthening our
exactness conditions, however, we can nevertheless derive a bound by comparing
to Simpson's rule:

\begin{lem}
Let $n=3$, $\xi_0=0$, $\xi_1=\frac{1}{2}$, $\xi_2=1$, and let $\hat{I}$ be
symmetric and affinely exact on some function $f^*\in C^4(\Rset,R)$ with
$o(h^4)$ errors for $f\in C^4(a,a+h)$. Then for $f\in C^4(a,a+h)$ we have
\begin{equation}
|\hat{I}[f]-I[f]|\le\frac{h^5}{2880}
\sup_{\xi\in[a;a+h]}|f^{(4)}(\xi)-\Omega_f|
\label{eq:errboundh5}
\end{equation}
where
\begin{equation}
\Omega_f = \frac{2880}{h^4}\left(\frac{f(a)+4f(a+\frac{h}{2})+f(a+h)}{6} -
q(f(a),f(a+\frac{h}{2}),f(a+h))\right)
\end{equation}
is bounded for $h\to 0$ because both $\hat{I}$ and Simpson's rule have errors of
order $o(h^4)$.
\begin{proof}
Defining
\begin{align*}
p_f(x) & = \sum_{i=0}^2 f(a+\xi_ih)L_i(x) + \Omega_f N(x) + M(x)\\
L_i(x) & = \frac{\prod_{j\not=i} (x-a-\xi_jh)}{h^2 \prod_{j\not=i}(\xi_i-\xi_j)} \\
N(x) & = \frac{1}{24}(x-a)(x-a-h)(x-a-\frac{h}{2})^2 \\
M(x) & = 4\frac{f(a+h)-f(a)-hf'(a+\frac{h}{2})}{h^3} (x-a)(x-a-h)(x-a-\frac{h}{2})
\end{align*}
we have
$e_f(a)=e_f(a+h)=e_f(a+\frac{h}{2})=e'_f(a+\frac{h}{2})=0$,
$\hat{I}[p_f]=I[p_f]$ by construction. Applying lemma \ref{lem:greenbound} with
$r=4$, $D=\frac{\rmd^4}{\rmd x^4}$, $H=\{1\}$, then yields $Df(\xi)=f^{(4)}(\xi)$,
$Dp_f(\xi)=\Omega_f$, and $C=\frac{h^5}{2880}$.
\end{proof}
\end{lem}

We note that for $\hat{I}$ equal to Simpson's rule, this is just the well-known
error bound on the latter, since $\Omega_f=0$ in this case.

Applying this result to the quadrature formula $q_2$ of eq.~(\ref{eq:q4def})
then completes the proof of Theorem \ref{thm:construct}.

Going beyond these easy, if perhaps not all that useful, results towards
tighter bounds that take the structure of the higher-order
nonlinear quadrature formula into account would likely require
additional constraints on $\hat{I}$ in order to ensure the existence of a
suitable family of functions that can be used to interpolate $f$ with a
function on which $\hat{I}$ is exact.

\section{Traditional quadrature rules as linear approximations}

We first note that when applying the construction of eq.~(\ref{eq:construct})
to the identity function $f^*(x)=x$ with inverse $f^*{}^{(-1)}(x)=x$ and
antiderivative $F^*(x)=\frac{1}{2}x^2+C$, we obtain the trapezoidal rule,
\begin{equation}
q(f_0,f_1) = \frac{\frac{1}{2} f_1^2-\frac{1}{2}f_0^2}{f_1-f_0} =
\frac{f_0+f_1}{2},
\end{equation}
which is affinely exact by construction, and hence is exact on all first-order
polynomials $f(x)=\alpha x+\beta$, whose second derivative vanishes identically.
We have thus given an alternative derivation of a well-known result:

\begin{corol}
The trapezoidal rule
\begin{equation}
\hat{I}[f] = \frac{h}{2}\left(f(a)+f(a+h)\right)
\end{equation}
is exact for all first-order polynomials $f(x)$ and satifies
\begin{equation}
\label{eq:errboundtrap}
|\hat{I}[f]-I[f]|\le\frac{h^3}{12}\lVert f''\rVert_\infty
\end{equation}
for $f\in C^2(\Rset)$.
\end{corol}

We note that when the trapezoidal rule is used for $q_1$ in
eq.~(\ref{eq:constructh4}), the rule $q_2$ obtained in this way is precisely
Simpson's rule, yielding an alternative proof of another well-known result:

\begin{corol}
The quadrature rule
\begin{equation}
\hat{I}[f] = \frac{h}{6}\left(f(a)+4f(a+\frac{h}{2})+f(a+h)\right)
\end{equation}
satisfies
\begin{equation}
\label{eq:errboundsimps}
|\hat{I}[f]-I[f]|\le\frac{h^5}{2880}\lVert f^{(4)}\rVert_\infty
\end{equation}
for $f\in C^4(\Rset)$.
\end{corol}

Similarly, the higher-order Newton-Cotes
rules can be obtained without any explicit reference to polynomial
interpolation by linearly combining the different evaluations of the integral
from $a$ to $a+h$ that can be formed using the trapezoidal rule on the nodes of
the higher-order Newton-Cotes rule and optimizing the coefficients of the
linear combination to minimize the total error:

\begin{corol}
Let $n$ be odd. Then the $n$-point quadrature rule
\begin{equation}
\hat{I}[f] = h\sum_{k=0}^{n-2} \alpha_k
\left(\xi_k\frac{f(a)+f(a+\xi_kh)}{2}+(1-\xi_k)\frac{f(a+\xi_kh)+f(a+h)}{2}\right)
\end{equation}
with $\xi_k=\frac{k}{n-1}$ and $\alpha_k$ given by the solution of the linear
equation system
\begin{equation}
\sum_{k=0}^{n-2}\alpha_k\left(\xi_k^j-\xi_k+1\right)=\frac{2}{j+1} \qquad
j=2,\ldots n
\end{equation}
is identical to the $n$-point Newton-Cotes rule with nodes $\xi_k$ (taking
$\xi_{n-1}=1$).
\begin{proof}
First, we note that by linearity, we can write
\begin{equation}
\hat{I}[f] = h\sum_{k=0}^{n-1} w_k f(a+\xi_k h)
\end{equation}
determining a unique linear $n$-point quadrature formula. Next, we note that
\begin{align}
\xi_kh\frac{f(a)+f(a+\xi_kh)}{2}+(1-\xi_k)h\frac{f(a+\xi_kh)+f(a+h)}{2}
&= \\ \nonumber f(a)h + \frac{1}{2}f'(a)h^2+\sum_{j=2}^n
\frac{1}{2j!}f^{(j)}(a)\left[\xi_k^j-\xi_k+1\right]h^{j+1}
+o(h^{n+1})
\end{align}
and hence demanding that the Taylor expansion of $\hat{I}[f]$ matches that of
$I[f]$ up to order $h^{n+1}$ amounts to the $(n-1)\times(n-1)$ linear equation
system
\begin{equation}
\sum_{k=0}^{n-2}\alpha_k\left(\xi_k^j-\xi_k+1\right)=\frac{2}{j+1} \qquad
j=2,\ldots n
\end{equation}
which has a unique solution since the rows are polynomials of different orders
in $\xi_k$ and hence must be linearly independent. This unique solution yields
an $n$-point quadrature formula that is exact on polynomials of order $n-1$
(whose derivatives from the $n$\textsuperscript{th} on all vanish), and hence
must be identical to the $n$-point Newton-Cotes rule on nodes $\xi_k$, which is
defined by this exactness.
\end{proof}
\end{corol}

One easily verifies that Simpson's ($n=3$), Boole's ($n=5$) and Weddle's
($n=7$) rules are recovered in this way.

\section{Explicit Non-linear Examples}

The construction of Theorem \ref{thm:construct} for $f^*(x)=\rme^x$ yields
the very quadrature rule whose empirical use by practitioners was the
motivation of the present inquiry, \emph{viz.}
\begin{equation}\label{eq:qexp}
q_1(\hat{f}) = \frac{f(a+h)-f(a)}{\log \frac{f(a+h)}{f(a)}}
\end{equation}
which by construction is symmetric and affinely exact for $f^*(x)=\rme^x$.
Since $\lambda\rme^{\alpha x}=\rme^{\alpha x+\beta}$ with $\beta=\log\lambda$,
the corresponding non-linear quadrature formula is also scalably exact
on all functions of the form $f(x)=\rme^{\alpha x}$. Finally, this quadrature
formula is quasilinear since the numerator is linear and any scalar factor
$\lambda$ cancels within the denominator.

In terms of errors, this rule satisfies
\begin{align}
|\hat{I}[f]-I[f]| &\le \frac{h^3}{12} \sup_{\xi\in[a;a+h]}
|f''(\xi)-\alpha_f^2\rme^{\alpha_f x+\beta_f}|\\
|\hat{I}[f]-I[f]| &\le
\left(\frac{h}{\alpha_f^2}-\frac{2}{\alpha_f^3}\tanh\frac{h\alpha_f}{2}\right) \sup_{\xi\in[a;a+h]}
|\alpha_f^2 f(\xi)-f''(\xi)|
\end{align}
since $f^*(x)=\rme^x$ satisfies $f(x)-f^*{}''(x)=0$, and
$C_f=\frac{h}{\alpha_f^2}-\frac{2}{\alpha_f^3}\tanh\frac{h\alpha_f}{2}$
by direct calculation. We note that
$C_f=\frac{h^3}{12}-\frac{\alpha_f^2h^5}{120}+O(h^7)$.

In applications like the Padé-Laplace method \cite{Yeramian:1987nature}
we also require the momenta of multiexponential functions, and thus need to 
integrate products of multiexponential functions given as data points and 
monomials $x^n$. In this
case, the integrand decays exponentially at large $x$, but grows polynomially
at small $x$, so that nonlinear quadrature rules for exponentials will only
work well at large $x$, while Newton-Cotes rules will be more appropriate at
small $x$. One way to determine where the change in regime to exponential decay
happens would be to consider numerical derivatives of the data $f_k=f(kh)$ and
to use the nonlinear quadrature rule only in the convex decaying region where
$f_k<f_{k-1}$ and $2f_k<f_{k-1}+f_{k+1}$, and to use Simpson's rule otherwise.

We can, however, do better than this by considering the integration-by-parts
identity
\begin{equation}
\int x^n \rme^{\alpha x}\rmd x = \frac{x^n}{\alpha}\rme^{\alpha
x}-\frac{n}{\alpha} \int x^{n-1}\rme^{\alpha x}\rmd x
\end{equation}
with solution
\begin{equation}
\int x^n \rme^{\alpha x}\rmd x = \sum_{k=0}^n \frac{(-1)^kn!}{(n-k)!}
\frac{x^{n-k}\rme{\alpha x}}{\alpha^{k+1}}
\end{equation}
and use the same heuristic that originally led us to consider
eq.~(\ref{eq:qexp}) in the first place to arrive at a quadrature formula for
moments of functions that are close to an exponential,
\begin{equation}
\hat{I}^n[f] = \sum_{k=0}^n
\frac{(-1)^kn!}{(n-k)!}\frac{(a+h)^{n-k}f(a+h)-a^{n-k}f(a)}{\log^{k+1}
\frac{f(a+h)}{f(a)}}h^{k+1}.
\end{equation}
This formula is scalably exact on functions of the form $f(x)=\rme^{\alpha x}$,
$\hat{I}^n[\lambda\rme^{\alpha x}]=I[x^n\lambda\rme^{\alpha x}]$, and has
$O(h^3)$ errors,
\begin{equation}
\left|\hat{I}^n[f]-\int_a^{a+h}x^nf(x)\right| = \frac{f(a)f''(a)-f'(a)^2}{12
f(a)}a^nh^3 + O(h^4),
\end{equation}
making it well-suited for usage with the Padé-Laplace method.

Finally, we note that for the case (common in applications) where we need to
estimate an improper integral out to infinity from a finite number of samples
$f(x_k)$, we can readily generalize eq.~(\ref{eq:qexp}) in order
to get a quadrature rule for improper integrals of the form
\begin{equation}
\int_a^\infty f(x)\,\rmd x \approx h\frac{f(a)}{\log \frac{f(a)}{f(a+h)}}
\end{equation}
assuming that $f$ is monotonically decaying with exponential speed such that
the integral converges and $f(a)>f(a+h)$. This rule is scalably exact on all
functions of the form $f(x)=\rme^{\alpha x}$, $\alpha<0$,
but the error analysis for the
case of the proper integral does not carry through (as there is notably no $h$
dependence of the left-hand side).

Given that the construction of Theorem \ref{thm:construct} takes the quotient
of two differences, the question of its numerical stability naturally arises.
This is even more pronounced in the case of the higher-order rule, which
involves an additional difference due to the negative coefficient
$-\frac{1}{3}$. 
As the examples of the trapezoidal rule and Simpson's rule in the previous
section show, numerical instability is not a given since there may be a
manifestly stable form of the quadrature rule that is mathematically equivalent
to the construction of Theorem \ref{thm:construct} in exact arithmetic.
Whenever possible, it is therefore desirable to bring the quadrature rule
derived from Theorem \ref{thm:construct} into a form that involves as few
differences as possible.

\section{Numerical Experiments}

We have tested the accuracy of eq.~(\ref{eq:qexp}), and of the corresponding
three-point quadrature formula by comparison with the trapezoidal rule and
Simpson's rule, respectively. To this end we consider the integrals of
\begin{align}
f_1(x) &= \rme^{-x}+\frac{1}{2}\rme^{-2x}, \\
f_2(x) &= \sum_{n=1}^\infty \rme^{-n x} = \frac{1}{\rme^x-1}, \\
f_3(x) &= \cosh x, \\
f_4(x) &= \sin x,
\end{align}
over the intervals $[a;b] = [0;1/2]$, $[3/2;2]$, $[0;1/2]$, and
$[\frac{\pi}{6};\frac{\pi}{3}]$, respectively.

\begin{figure}[p]
{\centering
\includegraphics[width=0.49\textwidth]{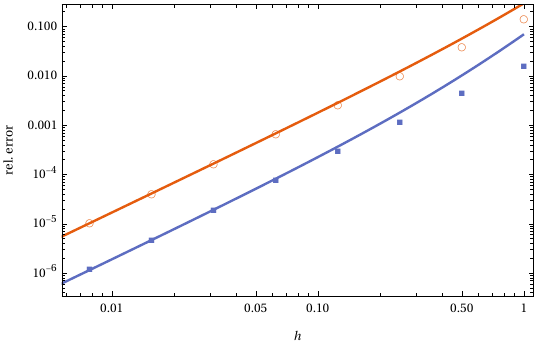}
\includegraphics[width=0.49\textwidth]{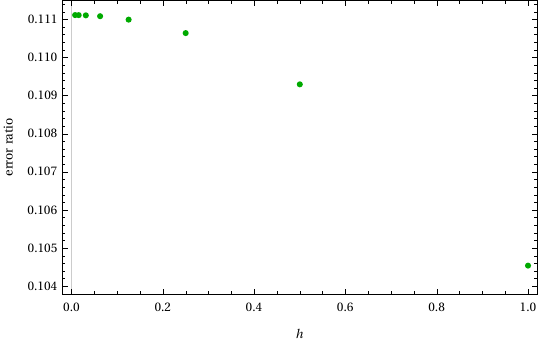}\\
\includegraphics[width=0.49\textwidth]{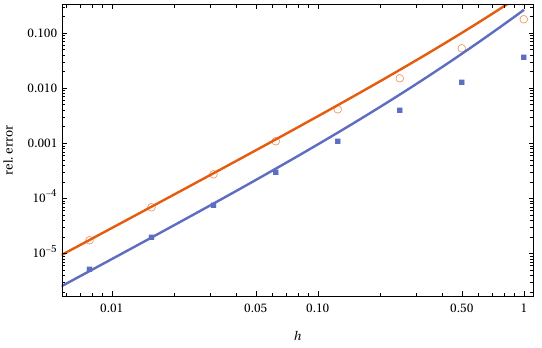}
\includegraphics[width=0.49\textwidth]{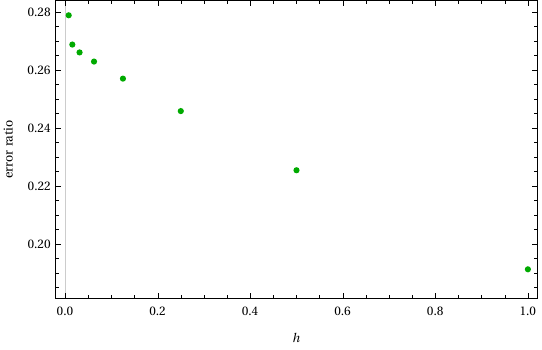}\\
\includegraphics[width=0.49\textwidth]{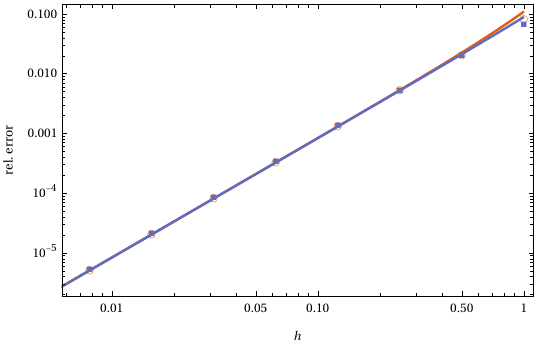}
\includegraphics[width=0.49\textwidth]{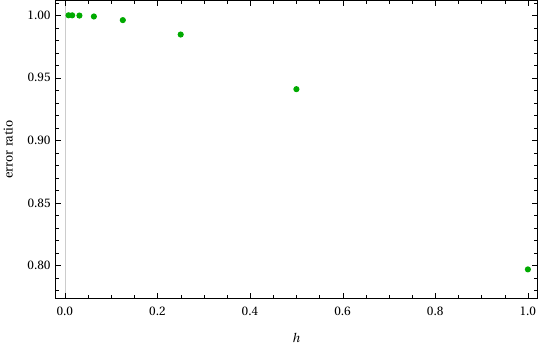}\\
\includegraphics[width=0.49\textwidth]{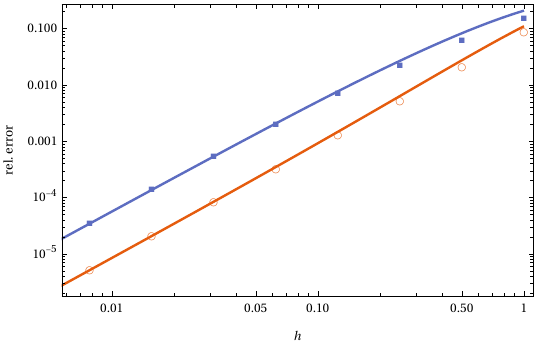}
\includegraphics[width=0.49\textwidth]{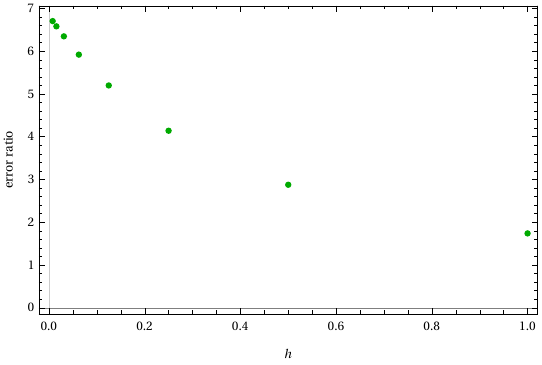}\\
}
\caption{\label{fig:trap}
Comparison between the non-linear exponential rule (solid blue squares) and the trapezoidal
rule (open red circles) of the relative error $|\hat{I}[f]-I[f]|/|I[f]|$ (left) and the
ratio of the errors between the two rules (right) on a range of integrands (top
to bottom: $\rme^{-x}+\frac{1}{2}\rme^{-2x}$, $[\rme^x-1]^{-1}$, $\cosh x$,
$\sin x$). See the text for details.
}
\end{figure}

\begin{figure}[p]
{\centering
\includegraphics[width=0.49\textwidth]{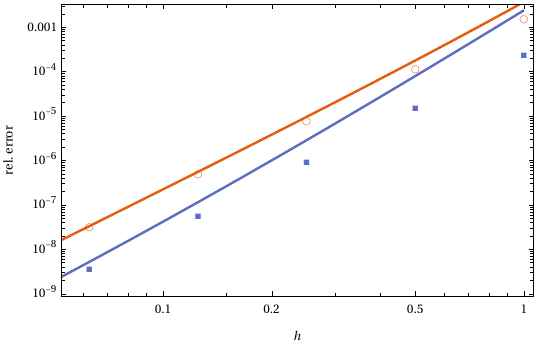}
\includegraphics[width=0.49\textwidth]{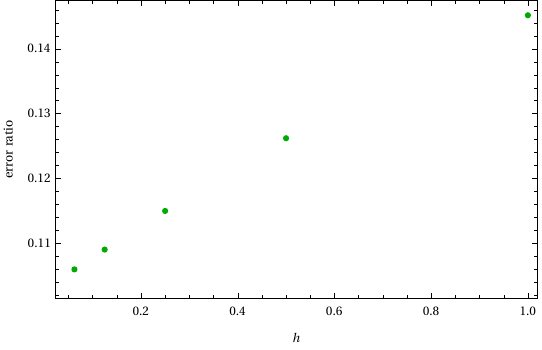}\\
\includegraphics[width=0.49\textwidth]{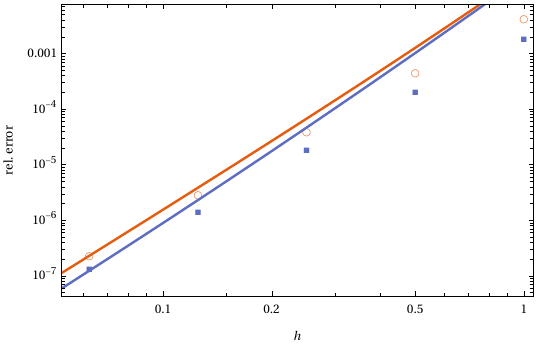}
\includegraphics[width=0.49\textwidth]{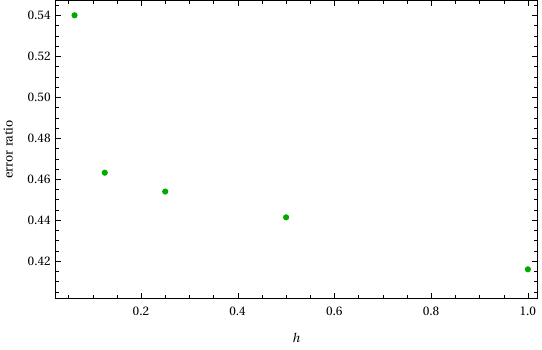}\\
\includegraphics[width=0.49\textwidth]{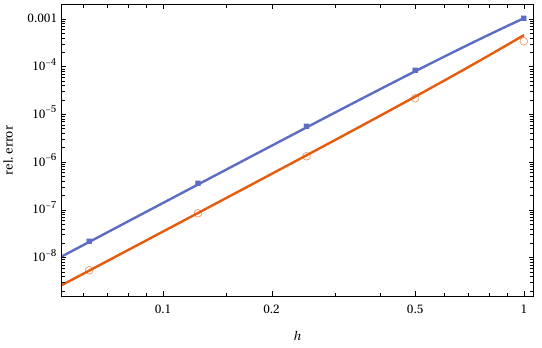}
\includegraphics[width=0.49\textwidth]{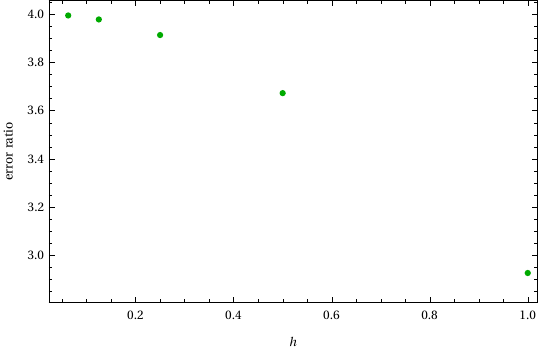}\\
\includegraphics[width=0.49\textwidth]{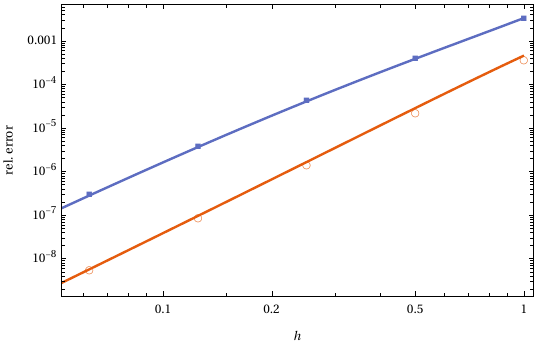}
\includegraphics[width=0.49\textwidth]{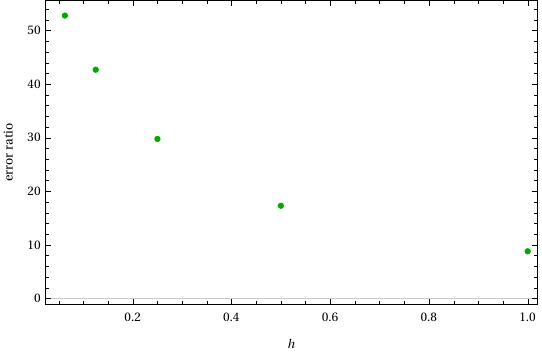}\\
}
\caption{\label{fig:simps}
Comparison between the higher-order non-linear exponential rule (solid blue
squares ) and
Simpson's rule (open red circles) of the relative error $|\hat{I}[f]-I[f]|/|I[f]|$ (left) and the
ratio of the errors between the two rules (right) on a range of integrands (top
to bottom: $\rme^{-x}+\frac{1}{2}\rme^{-2x}$, $[\rme^x-1]^{-1}$, $\cosh x$,
$\sin x$). See the text for details.
}
\end{figure}

First, we consider the accuracy of a single step $[a;a+h]$. Figures
\ref{fig:trap} and \ref{fig:simps} show the comparison between the relative
errors
\begin{equation}
e(h) = \frac{|\hat{I}[f]-I[f]|}{|I[f]|}
\end{equation}
for the $O(h^3)$ and $O(h^5)$ cases in their respective left columns. The right
columns show the error ratios
\begin{equation}
r(h) = \frac{|\hat{I}[f]-I[f]|}{|\hat{I}_0[f]-I[f]|}
\end{equation}
where $\hat{I}$ and $\hat{I}_0$ are the non-linear and linear quadrature
formulae, respectively.

It can be seen that on $f_1$ and $f_2$, which are well described as being
dominated by a leading exponential, the non-linear quadrature formulae
outperform their linear counterparts by approximately an order of magnitude in
error in the case of $f_1$ and a factor of $2$ to $4$ in the case of $f_2$. In the
case of $f_3$, which is a sum of two exponentials, but with opposite signs of
the exponent, the advantage of the non-linear rules is very small in the
$o(h^2)$ case and non-existent in the $o(h^4)$ case, where Simpson's rule is
more efficient by a factor of $4$. In the case of $f_4$, which not a sum of
real exponentials at all (although it is the sum of complex exponentials), the
non-linear rules tailored to real exponentials fare very poorly, with one to
two orders of magnitude larger errors than their linear counterparts (a large
part of this likely being due to the fact that $f_4$ is positive and concave on
the integration interval, while the non-linear rules implicitly assume a
positive convex, or negative concave, function).

Shown alongside the data in each case are the error bounds of
eqs.~(\ref{eq:errbound2}) and (\ref{eq:errboundh5}) for the non-linear rules
and eqs.~(\ref{eq:errboundtrap}) and (\ref{eq:errboundsimps}) for the
Newton-Cotes rules.

\begin{figure}[p]
{\centering
\includegraphics[width=0.49\textwidth]{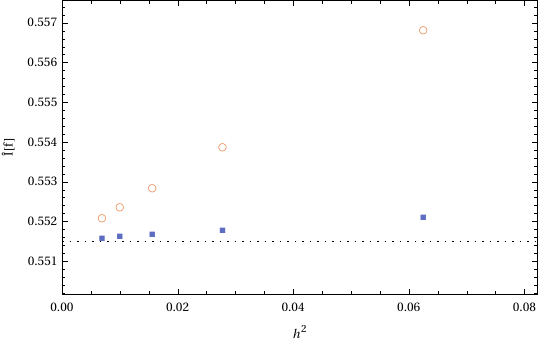}
\includegraphics[width=0.49\textwidth]{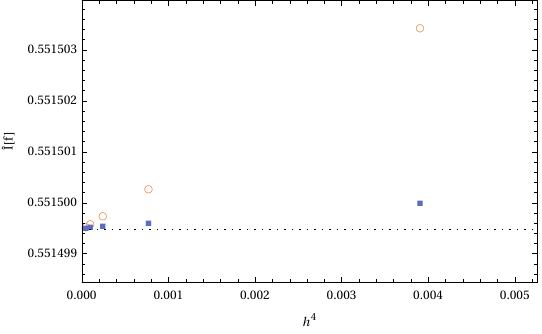}\\
\includegraphics[width=0.49\textwidth]{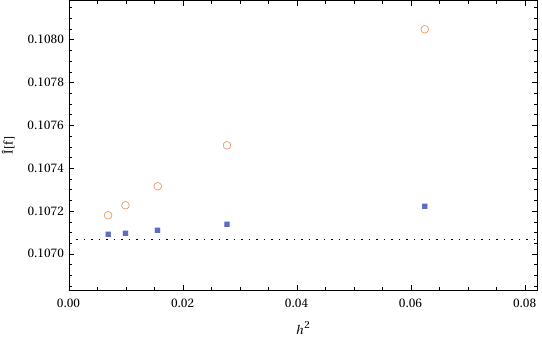}
\includegraphics[width=0.49\textwidth]{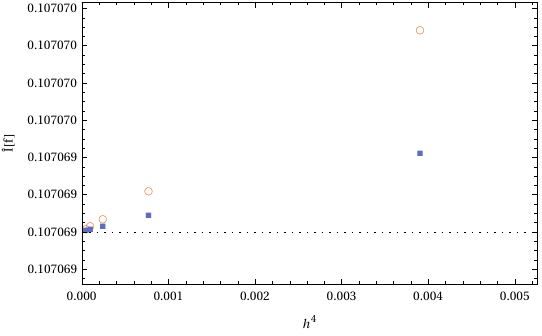}\\
\includegraphics[width=0.49\textwidth]{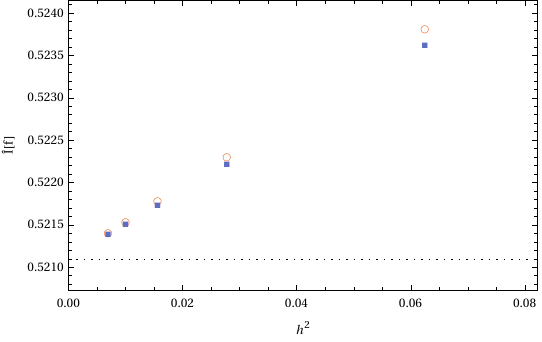}
\includegraphics[width=0.49\textwidth]{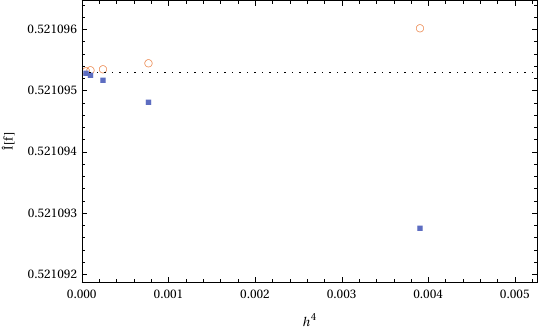}\\
\includegraphics[width=0.49\textwidth]{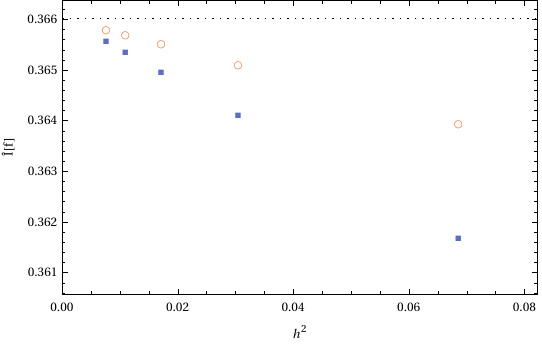}
\includegraphics[width=0.49\textwidth]{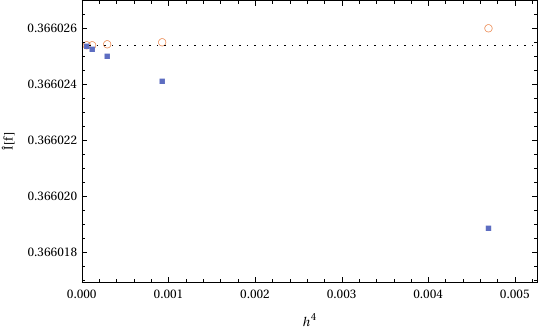}\\
}
\caption{\label{fig:multi}
Comparison between the results from multistep evaluation of integrals using the
non-linear exponential rule (solid blue squares) and the trapezoidal
rule (open red circles) on the left, and the higher-order non-linear
exponential rule (solid blue squares) and
Simpson's rule (open red circles) on the right. The integrands are (top
to bottom) $\rme^{-x}+\frac{1}{2}\rme^{-2x}$, $[\rme^x-1]^{-1}$, $\cosh x$,
$\sin x$. Note the differences in scale. See the text for details.
}
\end{figure}

Next, we consider the convergence of the full integrals over $[a;b]$ evaluated
using $N$ steps from $a+kh$ to $a+(k+1)h$ with $h=\frac{b-a}{N}$. Figure
\ref{fig:multi} shows the exact values of the integrals as dotted horizontal
lines and the numerical evaluations using the $O(h^2)$ and $O(h^4)$ multistep
rules in the left and right columns, respectively (note that the vertical
scales in the left and right columns differ markedly). The approach to the
continuum limit is as expected from the previous two figures.

\section*{Conflicts of interest}

The author is not aware of any personal, professional, financial, political or
other circumstances that could give rise to a relevant conflict of interest. 

\section*{Acknowledgements}

The author is grateful to Ron R. Horgan and Harvey B. Meyer for useful
comments, and wishes to thank Christian Remling for his MathOverflow answer
\cite{Remling:MO463802} pointing out the usefulness of the Green function in
obtaining bounds on the $L^1$ norm of the solution of a boundary value problem.

\end{document}